\documentclass[12pt, oneside]{amsart}   	
\usepackage{geometry}
\usepackage{amsmath, amssymb}
\usepackage{mathrsfs}
\usepackage{epsfig}
\usepackage{epic,eepic}              		
\usepackage{graphicx}				
\newtheorem{thm}{Theorem}[section]
\newtheorem{prop}[thm]{Proposition}
\newtheorem{prob}[thm]{Problem}
\newtheorem{ques}[thm]{Question}
\newtheorem{lem}[thm]{Lemma}

\newtheorem{conj}[thm]{Conjecture}
\newtheorem{rem}[thm]{Remark}

\theoremstyle{definition}

\numberwithin{equation}{section}

\title{Polynomials with dense zero sets and discrete models of the Kakeya conjecture and the Furstenberg set problem}
\author{Ruixiang Zhang}
\usepackage{fancyhdr}
\begin{document}

\maketitle
\begin{abstract}
We prove the discrete analogue of Kakeya conjecture over $\mathbb{R}^n$. This result suggests that a (hypothetically) low dimensional Kakeya set cannot be constructed directly from discrete configurations. We also prove a generalization which completely solves the discrete analogue of the Furstenberg set problem in all dimensions. The difference between our theorems and the (true) problems is only the (still difficult) issue of continuity since no transversality-at-incidences assumptions are imposed. The main tool of the proof is a theorem of Wongkew \cite{wongkew2003volumes} which states that a low degree polynomial cannot have its zero set being too dense inside the unit cube, coupled with Dvir-type polynomial arguments \cite{dvir2009size}. From the viewpoint of the proofs, we also state a conjecture that is stronger than and almost equivalent to the (lower) Minkowski version of the Kakeya conjecture and prove some results towards it. We also present our own version of the proof of the theorem in \cite{wongkew2003volumes}. Our proof shows that this theorem follows from a combination of properties of zero sets of polynomials and a general proposition about hypersurfaces which might be of independent interest. Finally, we discuss how to generalize Bourgain's conjecture to high dimensions, which is closely related to the theme here.
\end{abstract}
\section{introduction}

\subsection{Discrete models of Kakeya conjecture and the polynomial method}

A Kakeya set in $\mathbb{R}^n$ is a compact set which contains a unit line segment in every direction. The Kakeya conjecture concerning it, which has close connections to many harmonic analysis and combinatorial problems, states the following:
\begin{conj}[Kakeya Conjecture]\label{Kakeya}
A Kakeya set in $\mathbb{R}^n$ has full Hausdorff dimension.
\end{conj}

Conjecture \ref{Kakeya} is notoriously difficult, especially in high dimensions. Various partial results have been proved. One could ask an easier question about the Minkowski dimension instead and might obtain better lower bounds than the Hausdorff version. In either case, the best known lower bound for large $n$, which is due to Katz and Tao \cite{katz2002new}, is still far from $n$ (and yet highly nontrivial). These bounds are linear in terms of $n$ and significantly larger than easier bounds (which are around $\frac{n}{2}$) one could prove. For example, the best known Hausdorff dimension lower bound is $(2- \sqrt{2}) (n-4) + 3$ \cite{katz2002new}. In high dimensions, the arguments leading to these bounds exploit the additive structure of finite point sets in a line \cite{katz1999bounds} \cite{katz2002new}, which was in turn inspired by an earlier paper of Bourgain\cite{bourgain1999dimension}. In \cite{katz2002new}, the authors also proposed a combinatorial program (i.e. the ``$SD(\alpha)$'' conjecture there) that suggests one possible way to attack the full Kakeya.

To start doing quantitative analysis, we introduce some notations.When $A, B >0$, we use $A \gtrsim B$ to denote that $A \geq C B$ for the rest of the paper, where $C$ is an absolute constant. Similarly, if $\epsilon >0$ then $\gtrsim_{\varepsilon}$ means that a similar inequality holds where the implied constant depends only on $\varepsilon$, etc. The Minkowski version of the Kakeya conjecture, which asserts that the $\delta$-neighborhood of a fixed Kakeya set $K$ has volume $\gtrsim_{\varepsilon, K} \delta^{\varepsilon}$ for any $\varepsilon > 0$, can be viewed as a question about incidences between small balls and thin tubes in $\mathbb{R}^n$. Many authors have suggested to consider discrete models of this incidence problem. Perhaps the best known one was the Kakeya conjecture over finite fields, which was fully solved by Dvir\cite{dvir2009size} and now becomes a theorem.

\begin{thm}[Kakeya sets have full dimension in finite fields \cite{dvir2009size}]\label{Dvirsfinitetheorem}
For a finite field $\mathbb{F}$, if $K$ is a Kakeya set in $\mathbb{F}^n$, i.e. for each direction, there is a line $l \subseteq K$ parallel to that direction. Then $|K| \gtrsim_n |\mathbb{F}|^n$.
\end{thm}

Dvir's argument is now known as the ``polynomial method'' which has its roots in earlier number theory and combinatorics, and in turn inspired a lot of works, the best known one of which being the almost full solution to the Erd\"{o}s distinct distance conjecture in $\mathbb{R}^2$ by Guth and Katz\cite{guth2010erdos}. This method is not strong enough to settle the Kakeya problem in $\mathbb{R}^n$ due to subtle technical issues including the ``plany'' issue we will mention below. Nevertheless, it already shed some light on the real Kakeya conjecture: Guth\cite{guth2010endpoint} proved the endpoint multilinear Kakeya conjecture which was slightly better than the previous result by Bennett, Carbery and Tao\cite{bennett2006multilinear}. The most amazing part of Guth's work was that he was able to employ the similar but strengthened ``polynomial'' approach as in Dvir's paper to recover the previous result in \cite{bennett2006multilinear} which was proved in a totally different way. Moreover he got a slight improvement by allowing the endpoint case. The multilinear Kakeya was substantially weaker than Kakeya due to the ``transverse'' assumption. Therefore it does not quite suggest a way to the final solution of Kakeya. Yet it can already be of some use in the relevant harmonic analysis problems, see e.g. \cite{bourgain2011bounds}. Recently in an interesting paper \cite{guth2014restriction}, Guth gave an improved new restriction estimate via the polynomial method directly.

Back to the ``real'' Kakeya in $\mathbb{R}^n$, there is a natural formulation of a discrete model in $\mathbb{R}^n$, which will be one of the main topics we deal with in this paper. Let us first recall that the Kakeya conjecture can be viewed as an incidence problem between balls and tubes. It was perhaps Wolff \cite{wolff1999recent} \cite{wolff2003lectures} who first proposed a program to understand the analoguous problems of incidences between points and lines in $\mathbb{R}^n$ in order to understand at least some part of intuition towards the real Kakeya problem (in these references he also suggested the finite field Kakeya conjecture which was now Dvir's Theorem \ref{Dvirsfinitetheorem}). Indeed, it would at least seem strange to imagine that the worst case of incidences between tiny balls and thin tubes is much worse than the worst case of incidences between points and lines, and in vector spaces over finite fields the latter problem perfectly makes sense as well as in $\mathbb{R}^n$.

After Dvir's complete solution to the finite field Kakeya \cite{dvir2009size}, people started to go back to look at the relevant incidence problems in $\mathbb{R}^n$. In this direction, Wolff \cite{wolff1999recent} \cite{wolff2003lectures} already proposed a heuristic connection between Kakeya and a ``joints'' problem considered by various people (see e.g. \cite{chazelle1992counting} \cite{sharir1994joints} \cite{sharir2004point}). It has the extra ``transverse'' condition instead of the condition that the directions of the lines are evenly distributed as in the Kakeya problem (and has an accordingly changed conclusion). The philosophy of the connection is that a (hypothetical) Kakeya set of low dimension would have a lot of incidences between balls and tubes. If at many incidences there were $n$ ``transverse'' tubes then this would contradict the continuous version of the joints conjecture (as long as the latter could be formulated and proved). Indeed, Bennett, Carbery and Tao obtained a partial result conditioned on the angles of lines at the joint by their ``multilinear'' Kakeya estimate \cite{bennett2006multilinear}. But the story didn't end here: later after Dvir's work \cite{dvir2009size}, Guth and Katz \cite{guth2010algebraic} took up the polynomial method and proved the joints conjecture in dimension three.

The incidence between points and lines has a long history and a seminal result was the celebrated Szemer{\'e}di-Trotter theorem \cite{szemeredi1983extremal} which asserts that there are $\lesssim |P|^{\frac{2}{3}}|L|^{\frac{2}{3}} + |P| + |L|$ incidences between a finite set $P$ of points and a finite set $L$ of lines on the plane. There are examples suggesting that apart from the trivial bounds $|L|^2 + |P|$ and $|P|^2 + |L|$, this bound is essentially the best one we can hope. Thus we have a rather complete understanding of incidences between points and lines on the plane. In the higher dimensions, the picture is relatively vague and much work still has to be done. After the work of Dvir \cite{dvir2009size}, the polynomial method gradually became a powerful tool in this area. Indeed, there is a very simple proof of the Szemer{\'e}di-Trotter theorem by this method. Using the polynomial method, Guth and Katz \cite{guth2010algebraic} fully solved the joints problem in dimension three. Their arguments were simplified and generalized independently by Quilodr{\'a}n \cite{quilodran2010joints} and Kaplan-Sharir-Shustin \cite{kaplan2010lines} to arbitrarily high dimensions as the following theorem:

\begin{thm}[Joints Theorem \cite{quilodran2010joints} \cite{kaplan2010lines}]\label{jointstheorem}
For a set $L$ of lines in $\mathbb{R}^n$, a \emph{joint} is a point that lies on $n$ lines with linearly independent directions. Then the number of joints is $\lesssim_n |L|^{\frac{n}{n-1}}$.
\end{thm}

The next goal is to understand the ``plany'' issue of the Kakeya set: as we saw, much was known about the ``transversal'' incidences (note that in the multilinear Kakeya and the joints problem there are ``transversality'' assumptions built in), but it is more or less a common belief that the hardest part of the Kakeya set lie in the issue of ``plany'' incidences. In dimension two this issue does not exist and three would be the first interesting dimension. For the dimension three, Bourgain raised the following conjecture in AIM 2004 \cite{croot2004problems} as a discrete model of the real Kakeya problem. It is now a theorem also proved by Guth and Katz in their ``joints'' paper \cite{guth2010algebraic}:

\begin{thm}\label{Bourgainstheorem}
Assume that there are $N^2$ lines in $\mathbb{R}^3$ such that no $N$ lines lie on a common plane. If a point set $P$ has at least $N$ points lying on each of the previously given $N^2$ lines then $|P| \gtrsim N^3$.
\end{thm}

A naive attempt to generalize Bourgain's original conjecture to high dimensions is to ask the following:

\begin{ques}\label{naiveques}
Assume that there are $N^{n-1}$ lines in $\mathbb{R}^n$ such that no $N^{r-1}$ lines lie in a common $r$-dimensional affine linear subspace ($1 < r < n$). If a point set $P$ has at least $N$ points lying on each of the previously given $N^{n-1}$ lines then can we say $|P| \gtrsim_n N^n$?
\end{ques}

It may be striking to learn that, in fact, the answer is false in \emph{all sufficiently high dimensions}. Indeed, there are examples with all points and lines lying on a quadratic hypersurface such that they are somewhat more highly incidental than we expected. We will give a full account of this interesting story in our final section 5. In fact, similar phenomena for finite fields and for the continuous setting have been pointed out by Tao \cite{tao2005new}. We refer the readers there for these analogous phenomena. Also in other ranges than this ``Bourgain conjecture range (where the lines, points and incidences give data close to the situation in Question \ref{naiveques})'', highly incidental patterns on quadratic hypersurfaces happen even in dimension 4. See \cite{sharir2014incidence}\cite{solomon2014highly}.

We will elaborate this in the final section. There we will give a detailed counterexample that negatively answers Question \ref{naiveques} in sufficiently high dimensions. Our strategy would be taking a large number of the points and lines of low height. In the counting process, we use the classical sieve results \cite{birch1962forms} to prove the high number of incidences. We need a uniform bound for the number of points on any $r$-plane, where the uniformity is the key issue. We use an observation that goes back to Tao \cite{tao2005new} to reduce the counting into a quadratic one and then use a simple trick in a paper by Heath-Brown \cite{heath2002density}.

Based on the above concerns, one must be careful enough when trying to generalize Bourgain's conjecture (Theorem \ref{Bourgainstheorem}). The version we suggest is the following conjecture and problem:

\begin{conj}\label{highdimbourgainconj}
There exist integers $d_{n, 2}, \ldots, d_{n, n-1}$ such that: Assume that there are $N^{n-1}$ lines in $\mathbb{R}^n$ such that no $N^{r-1}$ lines lie in a common $r$-dimensional subvariety of degree $\leq d_{n, r}$ ($1 < r < n$). If a point set $P$ has at least $N$ points lying on each of the previously given $N^{n-1}$ lines then $|P| \gtrsim_n N^n$.
\end{conj}

\begin{prob}\label{highdimdegreeproblem}
If Conjecture \ref{highdimbourgainconj} holds, what is the minimal $d_{n, r}$?
\end{prob}

For example, Theorem \ref{Bourgainstheorem} implies Conjecture \ref{highdimbourgainconj} with $d_{3, 2} = 1 $ in Problem \ref{highdimdegreeproblem} in $\mathbb{R}^3$. In section 5 we will prove that in general $d_{n, n-1} \geq 2$ for large $n$.

Nevertheless, we should not forgot that this type of generalized questions have weaker assumptions than the analogues of Kakeya. We could not prove such a general conjecture in all dimensions (and the phenomena we will present in section 5 suggest that this conjecture might have very complicated aspects in high dimensions which are less related to the origional Kakeya conjecture). What we will do is to establish a partial result towards Conjecture \ref{highdimbourgainconj}, which has additional assumptions that make it a more natural discrete model of Kakeya conjecture:

\begin{thm}[Main theorem for the discrete model of Kakeya]\label{maththeoremofkakeya}
Given any fixed number $C, C'  > 0$. Assume that a set $L$ of lines in $\mathbb{R}^n$ satisfies that the direction (viewed as unit vectors) set of $L$ is either (a) a $\frac{C}{N}$-dense subset or (b) a $\frac{C}{N}$-separated subset with cardinality $|L| \geq C' N^{n-1}$ of the unit sphere $S^{n-1}$ (with respect to the Euclidean distance). If a point set $P$ has at least $N$ points lying on each line in $L$ then $|P| \gtrsim_{n , C, C'} N^n$.
\end{thm}

Here the direction set of a set of lines is the union of the directions of the lines. For simplicity we do not pass to projective spaces and adopt the convention that one line gives two opposite directions on the unit sphere. One can easily see that unlike the joints Theorem \ref{jointstheorem}, Theorem \ref{maththeoremofkakeya} does not require any local transverse condition and can indeed be viewed as a discrete version of the Kakeya problem in $\mathbb{R}^n$. As in Wolff's heuristic \cite{wolff1999recent} \cite{wolff2003lectures}, if one could replace points by small balls and replace lines by thin tubes (and allows an $\epsilon$ of loss on the exponent) then Theorem \ref{maththeoremofkakeya} will immediately become the statement that the Kakeya set has full Minkowski dimension. Unfortunately, we did not see any way of adapting the techniques of our proof to a solution to the Kakeya Conjecture \ref{Kakeya}. The main new difficulty in the real Kakeya would be the loss of ``precise'' incidence, leading to a complete loss of many algebraic structures of the polynomial we will construct.

Note that by passing to subsets we easily see that Conjecture \ref{highdimbourgainconj} implies Theorem \ref{maththeoremofkakeya}.

To conclude this part, we summarize the two main difficulties of the real Kakeya problem in $\mathbb{R}^n$: the ``plany'' issue and the ``continuous'' issue. The former issue is that most incidences might be plany which shuts down standard tools that exploit transversality and the latter issue is that the incidences of balls and tubes need not to be as accurate as the incidences of points and lines, which causes much trouble for algebraic treatments. For this issue, also the real Kakeya problem is significantly harder because that tubes can overlap a lot more than the lines, especially locally. Up to now, in high dimensions, the multilinear Kakeya estimates \cite{bennett2006multilinear}\cite{guth2010endpoint} can be established as a ``continuous'' version but always with transversality (at incidences) assumptions; the joints Theorem \ref{jointstheorem} has a weaker requirement on transversality but was only done in the discrete setting (recently in the three dimensional continuous setting there is a paper by Guth \cite{guth2014degree} in this direction); our Theorem \ref{maththeoremofkakeya} did not assume transversality at all but only achieves a discrete version. Despite these partial results, it is likely that substantially new ideas have to be invented in order to fully solve the real Kakeya conjecture in $\mathbb{R}^n$. Along the lines of this paper, later in this section we will propose a conjecture about polynomials that implies the Minkowski version of Kakeya. We will explain some possible ways to attack it but the solution almost surely requires new ingredients.

\subsection{A discrete model of the Furstenberg set problem}

Another interesting incidence problem in high dimensions is the discrete model of the Furstenberg set problem. The Furstenberg set problem is a generalization of Kakeya. It was likely to be inspired by Furstenberg's work \cite{furstenberg1970intersections} and was formulated by Wolff \cite{wolff1999recent} \cite{wolff2003lectures} and Tao \cite{taoedinburgh}. See these references and \cite{zhang2013Furstenberg} for an introduction. The problem is the following:

\begin{prob}[Furstenberg set problem]\label{furstenbergprob}
Fix $0 < \beta < 1$. If a compact set $S$ in $\mathbb{R}^n$ satisfies that for any direction $\omega \in S^{n-1}$, there is a line parallel to $\omega$ such that a $\beta$-dimensional subset of this line lies in $S$. Then what can we get as the best lower bound for the dimension of $S$?
\end{prob}

This problem can also be considered over finite fields and in both cases there are example sets with dimension $\leq \frac{n-1}{2} + \frac{n+1}{2} \beta$ \cite{zhang2013Furstenberg}. This is also my conjectural lower bound for $\mathbb{F}_p$ ($p$ prime) and $\mathbb{R}$, and for $\mathbb{R}$ this conjecture has a good reason (see below). However, provable lower bounds are not easy to obtain and in high dimensions we proved only a little bit of gain beyond the trivial bound for the ``critical exponent'' $\beta = \frac{1}{2}$ for $\mathbb{F}_p$\cite{zhang2013Furstenberg} (while for, say, $\mathbb{F}_{p^2}$, the trivial bound is not improvable for $\beta = \frac{1}{2}$, as Wolff already observed \cite{wolff1999recent}\cite{wolff2003lectures}).

One may also cook up a discrete model of Problem \ref{furstenbergprob} in $\mathbb{R}^n$. Interestingly enough, using the technique of this paper we can fully solve it up to constants. Namely, we will prove the following theorem:

\begin{thm}[Main theorem for the discrete model of Furstenberg]\label{maththeoremoffurstenberg}
Given any fixed number $C > 0$ and $0 \leq \beta \leq 1$. Assume that a set $L$ of lines in $\mathbb{R}^n$ satisfies that the direction set of $L$ is either (a) a $\frac{C}{N}$-dense subset or (b) a $\frac{C}{N}$-separated subset with cardinality $|L| \geq C' N^{n-1}$ of the unit sphere. If a point set $P$ has at least $N^{\beta}$ points lying on each line in $L$ then $|P| \gtrsim_{n , C} N^{\frac{n-1}{2} + \frac{n+1}{2} \beta}$. Moreover, the exponent $\frac{n-1}{2} + \frac{n+1}{2} \beta$ cannot be improved.
\end{thm}

Thus if we can ``replace'' the points and lines in Theorem \ref{maththeoremoffurstenberg} with balls and tubes (and allow an $\varepsilon$-loss on the exponent) we would obtain a lower dimension bound $\frac{n-1}{2} + \frac{n+1}{2} \beta$ for the Minkowski version of the Furstenberg set problem (and thus fully solve it). Hence following Wolff's heuristic, we put $\frac{n-1}{2} + \frac{n+1}{2} \beta$ as the conjectural lower dimension bound for the $n$ dimensional Furstenberg set problem. Note that in the two dimensional case, everything is concordant with Wolff's original work \cite{wolff1999recent} \cite{wolff2003lectures}. It is worth mentioning here that even in dimension 2 the Furstenberg set problem is widely open. The best known result there is a result of Bourgain \cite{bourgain2003erdos}, which is a nontrivial improvement of the trivial bound.

Obviously, Theorem \ref{maththeoremofkakeya} is a special case of Theorem \ref{maththeoremoffurstenberg}. Therefore we will only focus on the proof of Theorem \ref{maththeoremoffurstenberg} in the rest of this paper.

\subsection{Ingredients: Wingkew's theorem on polynomials with dense zero sets}

We come to the strategy of the proofs. Though Theorem \ref{maththeoremofkakeya} is a weaker version of the ultimate Conjecture \ref{highdimbourgainconj}, which is in turn a generalization of Theorem \ref{Bourgainstheorem}, our argument will be a ``global'' one unlike the ``local'' analysis (of singular points and planar points) in the proof of Theorem \ref{jointstheorem} and Theorem \ref{Bourgainstheorem} \cite{guth2010algebraic}\cite{kaplan2010lines}\cite{quilodran2010joints}. On the other hand, the basic strategy of our argument resembles Dvir's proof \cite{dvir2009size} more. Along the way of standard ``soft'' reasoning, we invoke a key theorem that was established in a totally different context \cite{wongkew2003volumes}.

\begin{thm}[A low degree polynomial cannot be too dense \cite{wongkew2003volumes}]\label{theoremondensepolynomial}
For $0 < \alpha < 1$ and positive integer $d$, the $\alpha$-neighborhood of the zero set of any nonzero polynomial of degree $d$ inside the unit cube $I^n \subseteq\mathbb{R}^n$ has volume $\lesssim_n \alpha d$.
\end{thm}

(In the original draft of this paper I was unaware of the references proving this theorem and gave my own version of the proof. After circulation of it I learned that this theorem was already known and there were at least two different proofs by Wongkew \cite{wongkew2003volumes} and Lotz \cite{lotz2012volume}. In particular this theorem was due to Wongkew. Nevertheless, our proof also has some new features that made us decide to keep it in the paper. See below.)

With Theorem \ref{theoremondensepolynomial}, we sketch the proof of Theorem \ref{maththeoremoffurstenberg} which will be done in the following section. First use the polynomial ham-sandwich theorem to partition the points evenly with a low degree polynomial $Q$. The highest homogeneous part $Q_h$ of $Q$ must vanish on the directions of all the lines that are contained in the zero set. But by Theorem \ref{theoremondensepolynomial} and a stereographic projection we know that the $\frac{\varepsilon}{\deg Q}$-neighborhood of those directions on the unit sphere can have arbitrarily small area when $\varepsilon$ is small. So by the assumptions, when $|P|$ is small there are a lot of lines that do not lie inside the zero set of $Q$. By a standard incidence counting argument, we bound the incidences between all the points and those lines and get a contradiction if there are too few points.

If one only wants to prove Theorem \ref{maththeoremofkakeya} which is the discrete analogue of Kakeya in $\mathbb{R}^n$, he does not need such a big machine explained above and the theorem follows directly from Theorem \ref{theoremondensepolynomial} and Dvir's argument \cite{dvir2009size}.

Along the lines of the proof of Theorem \ref{maththeoremoffurstenberg}, we propose a conjecture that is slightly stronger than the Minkowski version of Kakeya conjecture:

\begin{conj}\label{singulardirectionconj}
For any $\varepsilon > 0$, $H > 0$, a positive integer $N$ and a polynomial $Q$, we divide the unit cube $I^n \subseteq \mathbb{R}^n$ into cubes with edge length $N^{-1 - \varepsilon}$. For any direction $v \in S^{n-1}$, if there is a line that has the given direction and that there are $\geq \frac{1}{H} N^{1+ \varepsilon}$ small cubes satisfying: 1) each cube is ``roughly bisected'' (i.e. has the polynomial diving it into two parts, both of which have volume, say, $\geq \frac{1}{10}$ of the volume of the cube) by $Q$; 2) each cube intersects the line, then we say that the direction $v$ is \emph{$N, \varepsilon, H $-singular} with respect to this polynomial $Q$.

We conjecture that: Given any fixed numbers $\varepsilon, H > 0$. Then for any polynomial $Q$ whose degree is no more than a sufficiently large $N$, the set of $N, \varepsilon, H $-singular directions is not the whole $S^{n-1}$.
\end{conj}

The reason we present the conjecture in the current way is that we try to formulate a conjecture as weak as possible that implies Kakeya. In a later section, we will prove that this conjecture implies that Kakeya sets have Minkowski dimension $n$, with morally the same reasoning. In fact, we will see that by another trivial application of Theorem \ref{theoremondensepolynomial}, this conjecture is equivalent to a statement that is a very slight generalization of the Minkowski version of Kakeya!

Heuristically, the singular directions should be the directions ``near'' a point in the zero set of the highest homogeneous part of ``a small perturbation'' of $Q$ in $S^{n-1}$: Because there is a line with that direction such that the polynomial bisects ``too many'' cubes along the line, it is conceivable that it ``continues to bisecting cubes'' along it. However this is not very clear locally (i.e. for a single direction): We proved that a polynomial can ``turn very sharply'' \cite{zhang2013sharp}. Nevertheless, it's conceivable that ``few'' singular directions will behave ``strangely'' like that. If one manages to show that this is the case (which will require new techniques), Theorem \ref{theoremondensepolynomial} or any theorem of similar flavor might be useful to prove the conjecture.  In the above discussion we have to allow perturbation of $Q$ because a tiny high degree perturbation of the polynomial will behave quite like the original polynomial in the unit ball but behave arbitrarily at infinity.

After the formulation of Conjecture \ref{singulardirectionconj}, we give some evidence of it. For example we will prove a partial result in dimension 3 along the lines of some hairbrush-type ideas of Guth \cite{guth2014restriction} (hairbrush was introduced initially by Wolff\cite{wolff1995improved}). Not very surprisingly, we will see that this together with the degree reduction argument (for cubes) developed by Guth \cite{guth2014degree} can recover Wolff's $\frac{5}{2}$ bound (for Minkowski dimension). Basically a result better than our partial result will yield new results for Minkowski dimension of Kakeya sets in dimension 3. We leave this very challenging direction to the interested readers.

Next we sketch the proof of Theorem \ref{theoremondensepolynomial}. In \cite{wongkew2003volumes} and \cite{lotz2012volume} there are two different proofs of it and we will give our own version in section 4. Lotz's proof uses Weyl's tube formula and Crofton's theorem in integral geometry to calculate the volume of the neighborhood (and use the algebraicity to bound degrees of Gauss maps). Wongkew's proof and our proof have essentially the same input and are closely related. They are heuristically as follows: By the directed area argument in \cite{guth2010endpoint}, we know the area of the zero set inside $D^n$ is $\lesssim_n d$. One might then think that its $\alpha$-neighborhood has volume $\lesssim \alpha d$. However, the zero set might be very ``degenerate'' and behaves like ``low-dimensional'' objects and thus be very dense while having small surface area. To overcome this serious issue, both proofs use an induction on dimension and use a version of generalized Harnack curve theorem to prevent bad cases. The references of Harnack's curve theorem is \cite{harnack1876vieltheiligkeit} and what we will use as a generalization is the Oleinik-Petrovskii-Thom-Milnor version \cite{oleinik1949topology}\cite{thom1965homologie}\cite{milnor1964betti}. Wongkew set up a ``sufficiently dense grid of hyperplanes'' to achieve the conclusion while we choose a ``sufficiently large continuum of hyperplanes'' to do the proof. In our proof, we first prove a very general ``local'' lemma (Proposition \ref{inductionofdimprop}) for an arbitrary hypersurface and then make the theorem an easy corollary. We keep our proof in this paper because we thought that Proposition \ref{inductionofdimprop} might be useful in related problems.

We also notice that Theorem \ref{theoremondensepolynomial} is the best possible one can hope, since there is a constant $C_n$ such that for every large $d$ there exists a degree $d$ polynomial such that its zero set is $\frac{C_n}{d}$-dense inside the unit cube. In fact, it suffices to take a union of translated coordinate hyperplanes (which form a ``grid'').

\section*{Acknowledgements}
I was supported by Princeton University and the Institute for Pure and Applied Mathematics (IPAM) during the research. Part of this research was performed while I was visiting IPAM, which is supported by the National Science Foundation. I thank IPAM for their warm hospitality. I would like to thank Zeev Dvir, Jordan Ellenberg, Larry Guth, Roger Heath-Brown, Marina Iliopoulou, Nets Katz, Noam Solomon, Terence Tao and Ruobing Zhang for very helpful discussions. Zeev Dvir pointed out that we have to allow perturbation in the above remark after Conjecture \ref{singulardirectionconj}. The results in the final section was mainly inspired by the discussion with Terence Tao and Roger Heath-Brown. Terence Tao pointed out to me his work \cite{tao2005new}. Roger Heath-Brown gave me useful suggestions on the proof of Theorem \ref{counterexamplethm}. Marina Iliopoulou helped me to correct a mistake in the initial version. It was Larry Guth who brought \cite{lotz2012volume} to my attention.

\section{Proof of Theorem \ref{maththeoremofkakeya} and Theorem \ref{maththeoremoffurstenberg}}

As we pointed out above, Theorem \ref{maththeoremofkakeya} is a special case of Theorem \ref{maththeoremoffurstenberg}. Thus we will only prove Theorem \ref{maththeoremoffurstenberg}. We assume Theorem \ref{theoremondensepolynomial} in this section. First we recall the familiar polynomial partition lemma which can be proved in a standard way by the polynomial ham-sandwich theorem (see e.g. \cite{guth2010erdos})

\begin{lem}[Polynomial partition lemma]\label{polypartitionlemma}
Let $P$ be a finite set of points in $\mathbb{R}^n$, and let $D$ be a positive integer. Then there exists a nonzero polynomial $Q$ of degree $\leq D$ and a decomposition
\begin{equation}
\mathbb{R}^{n}=\{Q=0\}\sqcup U_{1} \sqcup \cdots \sqcup U_{m}
\end{equation}
Here $U_{1},\ldots ,U_{m}$ are open sets (which will be called \emph{cells}) bounded by $\{Q=0\}$, such that $m \sim_n D^n$ and that each cell $U_{i}$ contains $O_{n}(\frac{|P|}{D^n})$ points.
\end{lem}

\begin{proof}[Proof of Theorem \ref{maththeoremoffurstenberg}]
We only concentrate on case (a). Case (b) is similar. First we prove the lower bound. Using the polynomial partition Lemma \ref{polypartitionlemma}, we find a nonzero polynomial $Q$ of degree $\leq D$, where $D$ is a number to be determined later, such that the $\mathbb{R}^n$ is decomposed into its zero set $\{Q=0\}$ and cells $U_1, \ldots, U_m$ where $m \sim_n D^n$ and $|P \cap U_i| \lesssim_{n} \frac{|P|}{D^n}$.

We decompose the line set $L$ into two subsets $L_{alg}$ and $L'$ such that $L_{alg} = \{l \in L: Q \text{ vanishes on } l \}$ and $L' = L \setminus L_{alg}$. If a line $l \in L_{alg}$, then its direction gives two points $v_{l}$ and $-v_l$ on $S^{n-1}$. Note that the highest homogeneous part $Q_h$ of $Q$ must vanish on $v_{l}$ because of the following reason: $Q = Q_h + \text{ lower terms}$. Thus $Q(x + tv_l) (x \in \mathbb{R}^n, t \in \mathbb{R}) = t^{\deg Q} Q_h (v_l) + \text { lower terms in terms of  } t$. Since it vanishes identically (as a polynomial of $t$) for any $x \in l$, we deduce $Q_h (v_l) = 0$. Through the standard stereographic projection map $\phi$, there is a bi-Lipschitz diffeomorphism (which is also the restriction of an rational, in fact quadratic, morphism of bounded degree both on the numerator and on the denominator) between an absolute (i.e. independent of anything except the dimension $n$) open subset $U^{n-1} \subseteq S^{n-1}$ with Lipschitz boundary and the unit cube $I^{n-1} \subseteq \mathbb{R}^{n-1}$. Moreover by further restriction on $U^{n-1}$ we may assume this diffeomorphism has Jacobian $\sim_n 1$ uniformly. By abuse of notation we also call this restriction $\phi$. Now since $\phi$ has bounded degree on numerators and denominators, the zero set of ${\phi^{-1}}^* (Q_h)$ is the same as the zero set of a polynomial of degree $\lesssim_{n} D$. By Theorem \ref{theoremondensepolynomial}, we can choose $C_{n-1}$ such that for any $J > 0$, the $\frac{1}{C_{n-1} J \deg Q}$-neighborhood of the zero set of ${\phi^{-1}}^* (Q_h)$ inside the unit cube $I^{n-1}$ has ($n-1$-dimensional ) volume $\leq \frac{1}{J}$. Because the restriction $\phi|_{U^{n-1}}$ is bi-Lipschitz with Jacobian $\sim_n 1$ uniformly, we can start from a sufficiently large $J$ depending only on $n$ such that the following holds: there exists some $C_{n}' >0$ depending only on $n$ such that the $\frac{1}{C_{n}' \deg Q}$-neighborhood of the zero set of $Q_h$ inside $U^{n-1}$ has area $\leq \frac{1}{2} Area(U^{n-1})$. We will require that
\begin{equation}\label{reqofQ1}
\frac{1}{C_{n}' \deg Q} \geq \frac{2C}{N}
\end{equation}

Then by a simple area counting argument and the $\frac{C}{N}$- dense property of the directions of the lines in $L$, we easily deduce that (note that $U^{n-1}$ has nice Lipschitz boundary which has area $0$) when $N$ is sufficiently large (depending on $n$) there are $\gtrsim_{n} N^{n-1}$ lines in $L$, each of whose directions gives rise to a point in $U^{n-1} \setminus \{Q_h = 0\}$. i.e. $|L'| \gtrsim_n N^{n-1}$.

Since $\deg Q \leq D$, to satisfy (\ref{reqofQ1}) it suffices that $D \lesssim_n N$. From now on we will require
\begin{equation}\label{reqofD}
D \lesssim_n N
\end{equation}
such that the implied constant makes (\ref{reqofQ1}) hold. By the discussion of the previous paragraph, $|L'| \gtrsim_n N^{n-1}$. By discarding extra lines we may assume $|L'| \sim_n N^{n-1}$ (as in the following argument we only use the facts that $L' \ \subseteq L$ and that $L' \cap L_{alg} = \emptyset$).

We use $\mathcal{I} (P, L')$ to denote the number of \emph{incidences} between $P$ and $L'$. On one hand, we have $|L'| \gtrsim_n N^{n-1}$ and that each $l \in L'$ is incident to $\geq N^{\beta}$ points in $P$. Thus
\begin{equation}\label{incidence1}
\mathcal{I} (P, L') \gtrsim_{n} N^{n + \beta -1}.
\end{equation}

On the other hand, since $Q$ is of degree $\leq D$, each $l \in L'$ intersects $\{Q = 0\}$ for $\leq D$ times. Thus $l$ can only enter $\leq D$ cells. Look at all the incidences $(p, l) \in P \times L' $ such that $Q(p) = 0$ or that there does not exist another $p' \in P \cap l$ which lies in the same cell as $p$ does. By the fact that each line can enter $\leq D$ cells we see that the number of this type of incidences is $\lesssim |L'| D \lesssim_n N^{n-1} D $. We call all the other incidences ``type II''. Take a type II incidence $(p, l) \in P \times L'$. Assume that $l$ passes through $t\geq 2$ points in the cell where $p$ lies, then $l$ has $t$ type II incidences inside this cell while there are ${t \choose 2} \geq t$ pairs of different points that lie in the intersection of $l$ and this cell. Moreover, once a pair of different points are given, we can only have a unique line that passes through these two points. Therefore, the total number of type II incidences is no more than the total number of point pairs inside some same cell, which is $\lesssim_n m(\frac{|P|}{D^n})^2 \lesssim_n \frac{|P|^2}{D^n}$. Summing over the two types, we deduce that:
\begin{equation}\label{incidence2}
\mathcal{I} (P, L') \lesssim_{n} N^{n-1} D + \frac{|P|^2}{D^n}.
\end{equation}

Now we can take a very tiny constant $\delta_n$ depending only on $n$ and take $D = \lfloor\delta_n N^{\beta}\rfloor$ such that (\ref{reqofD}) always holds for all large $N$ and all $\beta \in [0, 1]$, and that when one compares (\ref{incidence1}) and (\ref{incidence2}) he obtains
\begin{equation}
N^{n + \beta -1} \lesssim_{n} \frac{|P|^2}{N^{n\beta}}.
\end{equation}

Therefore
\begin{equation}
|P| \gtrsim_n N^{\frac{n-1}{2} + \frac{n+1}{2} \beta}.
\end{equation}

Finally we show that the exponent $\frac{n-1}{2} + \frac{n+1}{2} \beta$ cannot be improved. We will construct an example similar to our arguments in \cite{zhang2013Furstenberg}, which in turn has its two dimensional root in Wolff's argument \cite{wolff1999recent} \cite{wolff2003lectures}. Without loss of generality, we assume $n \geq 2$ and that $N$ is sufficiently large.

Take a number $M$ to be determined later. Consider the point set
\begin{equation}
P_1 = \{ (\frac{a}{a+bM}, \frac{am_{1,1} + bm_{1,2}}{(a+bM)M}, \ldots, \frac{am_{n-1,1} + bm_{n-1,2}}{(a+bM)M}) : 1 \leq a, b \leq M^{\beta}, 1 \leq m_{i, j} \leq M\}
\end{equation}
 and the line set $L_1 = \{l_{m_{1,1}, m_{1,2}, m_{2,1}, m_{2,2}, \ldots, m_{n-1, 1}, m_{n-1, 2}}: 1 \leq m_{i, j} \leq M\}$ where the equation of the line $l_{m_{1,1}, m_{1,2}, m_{2,1}, m_{2,2}, \ldots, m_{n-1, 1}, m_{n-1, 2}}$ is
 \begin{equation}
x_2 = \frac{Mm_{1, 1} x_1 +m_{1, 2} (1-x_1)}{M^2}, x_3 = \frac{Mm_{2, 1} x_1 +m_{2, 2} (1-x_1)}{M^2}, \ldots, x_n = \frac{Mm_{n-1, 1} x_1 +m_{n-1, 2} (1-x_1)}{M^2}.
 \end{equation}

 Then there are $\geq M^{2\beta}$ points of $P_1$ on each line of $L_1$ (according to different possibilities of $a, b$). Also, the direction of the line $l_{m_{1,1}, m_{1,2}, m_{2,1}, m_{2,2}, \ldots, m_{n-1, 1}, m_{n-1, 2}}$  is parallel to the vector $(1, \frac{Mm_{1, 1} - m_{1, 2}}{M^2}, \frac{Mm_{2, 1} - m_{2, 2}}{M^2}, \ldots, \frac{Mm_{n-1, 1} - m_{n-1, 2}}{M^2})$. Note that the vectors $\{(1, \frac{Mm_{1} - m_{2}}{M^2}): 1 \leq m_1, m_2 \leq M\}$, when renormalized onto the unit circle, form a $\frac{C_1}{M^2}$-dense subset of an absolute open subset of the unit circle. Here $C_1$ is an absolute constant. Therefore the (renormalized) directions of the lines in the set $L_1$ form a $\frac{C_n}{M^2}$-dense subset of an absolute open subset of $S^{n-1}$. Here $C_n$ is a constant depending only on $n$. Hence the union of a finite number of rotated copies of $L_1$ (together with the according union of rotated $P_1$) will satisfy that the directions of the lines in it form a $\frac{C_n}{M^2}$-dense subset of $S^{n-1}$. Here the number of copies depends only on $n$. Thus if we choose $M = \lfloor HN^{\frac{1}{2}}\rfloor$ where $H$ is a sufficiently large constant depending only on $n, C$, then this union satisfies all the assumptions of Theorem \ref{maththeoremoffurstenberg}.

 We only need to count the total points in this configuration. The number of copies of $P_1$ depends only on $n$. Within $P_1$, there are $\lesssim_n (\max\{a_0, b_0\} M)^{n-1}$ points with fixed $a = a_0, b= b_0$. Taking the summation over $a$ and $b$, we deduce that the number of points in $P_1$ is $\lesssim_n M^{n-1 + (n+1) \beta}$. Since $M = \lfloor HN^{\frac{1}{2}}\rfloor$, the total number of points in this configuration is $\lesssim_{n, C} N^{\frac{n-1}{2} + \frac{n+1}{2} \beta}$. Therefore, the example shows that the exponent $\frac{n-1}{2} + \frac{n+1}{2} \beta$ cannot be improved.
\end{proof}

\begin{rem}
If one only wants to prove Theorem \ref{maththeoremofkakeya}, he does not need the polynomial partition Lemma \ref{polypartitionlemma}. In fact, a polynomial vanishing on the entire $P$ with lowest possible degree, together with Dvir's argument \cite{dvir2009size} and Theorem \ref{theoremondensepolynomial}, will work for him and produce a rather simple proof. Also the existence of an example configuration with $\sim_{n, C} N^n$ points is trivial.
\end{rem}

\section{Conjecture \ref{singulardirectionconj} is almost equivalent to the Minkowski version of Kakeya}

In this section, we prove that Conjecture \ref{singulardirectionconj} is slightly stronger than the statement that a Kakeya set must have full Minkowski dimension. First recall the polynomial ham-sandwich theorem mentioned in e.g. \cite{guth2010endpoint} which was originally proved in \cite{stone1942generalized}:
\begin{thm}[Polynomial ham-sandwich theorem]\label{Polyhamsandwich}
Let $M = {{n+d} \choose n} -1$. Let $U_1, \ldots, U_M$ be finite volume open sets in $\mathbb{R}^n$. There there is a nonzero real polynomial $Q$ of degree $\leq d$ whose zero set bisects each $U_i$. More precisely,
\begin{equation}
Vol (\{x \in U_i : Q(x) > 0\}) = Vol (\{x \in U_i : Q(x) < 0\}).
\end{equation}
\end{thm}

\begin{proof}[Proof that Conjecture \ref{singulardirectionconj} implies that Kakeya sets have full Minkowski dimension]
Any given Kakeya set $K$ is a compact set. We rescale the set to a compact set $K'$ inside the unit cube satisfying that there is a line segment of some fixed length $\delta = \delta_K$ parallel to every direction $\in S^{n-1}$ that is contained in $K'$. Fix any $\varepsilon > 0$. Take $H \gtrsim_n \frac{1}{\delta}$ and $N$ sufficiently large (depending on $\varepsilon, \delta, K$). It is easy to see that there is a constant $C_n$ such that all the small cubes of size $N^{-1-\varepsilon}$ that intersects $K'$ is contained inside the $C_n N^{-1- \varepsilon}$-neighborhood of $K'$. If the number of such cubes is $\lesssim_n N^n$ for a carefully chosen implied constant, then by the polynomial ham-sandwich Theorem \ref{Polyhamsandwich} there is a nonzero polynomial of degree $\leq N$ that bisects every such cube. But when $N$ is sufficiently large this obviously violates Conjecture \ref{singulardirectionconj}. Therefore the $C_n N^{-1- \varepsilon}$-neighborhood of $K'$ has volume $\gtrsim_n N^{-n\varepsilon}$ for all sufficiently large $N$. Since $\varepsilon$ is arbitrary, $K'$ has full Minkowski dimension and so does $K$.
\end{proof}

In fact by an application of Theorem \ref{theoremondensepolynomial}, we can prove that the reverse direction is ``almost'' true. We will prove the following:

\begin{thm}\label{eqthm}
Conjecture \ref{singulardirectionconj} is equivalent to the following statement: For any given $\delta > 0, \varepsilon > 0$, assume $I^n$ is divided into $N^n$ small cubes of side length $\frac{1}{N}$. If a set $P$ is a union of some cubes such that for each direction $v$, there exists a line $l_v$ in the direction such that the length measure of $P \bigcap l_v$ is $\geq \delta$ then the volume $|P| \gtrsim_{n, \delta, \varepsilon} N^{-\varepsilon}$.
\end{thm}

\begin{proof}
See the $N$ in the statement as ``$N^{1+\varepsilon}$'' in the conjecture. Then the above proof works to prove that the conjecture implies the statement. For the reverse direction, assuming every direction is singular, we deduce that the union of cubes that are roughly bisected by the polynomial satisfies the assumption of the statement for $\delta = \frac{1}{H}$ (and ``$N$'' replaced by $N^{1+ \varepsilon}$). Thus the $N^{-1-\varepsilon}$-neighborhood of the zero set of the polynomial inside $I^n$ must have a volume $\gtrsim_{n, H , \varepsilon} N^{- \frac{\varepsilon}{2}}$. By Theorem \ref{theoremondensepolynomial} the polynomial cannot have such a low degree as $N$. A contradiction.
\end{proof}

By standard arguments, a modification of the statement in Theorem \ref{eqthm}, with a slightly stronger requirement that the chosen cubes lining up on the line in each direction are consecutive, is equivalent to the (lower) Minkowski version of Kakeya.

To give some evidences that Conjecture \ref{singulardirectionconj} is reasonable, we prove it in the case that the zero set of $Q$ consists of $N$ hyperplanes. In fact, assuming that $N$ is large, then for each hyperplane $h_i$ with normal vector $n_i$ and any direction $v \in S^{n-1}$, any line of direction $v$ intersects $\lesssim_n \max \{\frac{1}{\frac{\pi}{2} - <n_i, v>}, N^{1+ \varepsilon}\}$ small cubes bisected by $h_i$ where $<n_i, v>$ is the angle between $\pm n_i$ and $ \pm v$. If every direction is $N, \varepsilon, H$-singular, then for every $v$ we have
\begin{equation}
\sum_{i=1}^N \max \{\frac{1}{\frac{\pi}{2} - <n_i, v>}, N^{1+ \varepsilon}\} \gtrsim_{n, H} N^{1 + \varepsilon}.
\end{equation}

Integrate this inequality over $v \in S^{n-1}$ and note that
\begin{equation}
\int_{S^{n-1}}  \max \{\frac{1}{\frac{\pi}{2} - <n_i, v>}, N^{1+ \varepsilon}\} \mathrm{d} Area S^{n-1} \lesssim_n \log N,
\end{equation}
we reach a contradiction for all large $N$.

However, this method already does not work for a union of quadratic hypersurfaces and it seems that more global properties have to be exploited.

The first interesting case for Conjecture \ref{singulardirectionconj} is dimension three. Now the ideas of Guth \cite{guth2014restriction} (which is Wolff's hairbrush method \cite{wolff1995improved}) can already make progress in this case. Essentially, the proof of Lemma 4.9 in \cite{guth2014restriction} implies:

\begin{thm}\label{3dprogress}
In $\mathbb{R}^3$, assume the unit cube is cut into $N^3$ small cubes of side length $\frac{1}{N}$. Assume a polynomial $Q$ has degree $D$, and that a set $V$ of $\frac{1}{N}$-separated directions such that for each direction $v \in V$, there is a line in that direction that passes through $\gtrsim N$ small cubes in which $Q$ has a zero.Then $|V| \lesssim D^2 N \log^2 N$.
\end{thm}

It is easy to see (and may not be surprising) that Theorem \ref{3dprogress} together with the degree reduction arguments in \cite{guth2014degree} can recover Wolff's $\frac{5}{2}$ dimension bound.

If we can replace the $D^2 N \log^2 N$ by $D^{1+ \varepsilon} N$ in Theorem \ref{3dprogress} then the dimension 3 case of Conjecture \ref{singulardirectionconj} and thus the Minkowski version of Kakeya conjecture can be proved. Also any improvement where the power of $D$ is smaller than $2$ would lead to nontrivial progress of the Kakeya problem in dimension three.

\section{A new proof of Theorem \ref{theoremondensepolynomial}}
In this section, we give our own proof of Theorem \ref{theoremondensepolynomial} which, as we explained, asserts that a low degree polynomial can't have its zero set being too dense inside the unit cube. Two other proofs can be found in \cite{wongkew2003volumes} and recently \cite{lotz2012volume}. The strategy of our proof was already explained in the introduction. We will need some preparation before the proof.

We will perform an induction on dimension and thus we sometimes also use $k$ to denote an arbitrary dimension (in applications it will be between $1$ and $n$). First recall a version of Harnack's curve theorem in high dimensions followed by the results of Oleinik-Petrovskii-Thom-Milnor \cite{oleinik1949topology}\cite{thom1965homologie}\cite{milnor1964betti}:

\begin{thm}[Harnack's curve theorem for hypersurfaces in high dimension]\label{generalizedharnack}
The zero set of a polynomial in $\mathbb{R}^k$ with degree $d$ has $\lesssim_k d^k$ connected components.
\end{thm}

Next we recall the concept of \emph{directed area}. For this part, the reader can see \cite{guth2010endpoint} (there it was called \emph{directed volume} and here we change the name slightly to emphasize that it is a hypersurface measure). For the sake of completeness, we will review everything we need here. Note that the singular set of any real algebraic hypersurface in $\mathbb{R}^k$ has zero hypersurface measure (one can make this rigorous by taking, say, the $k-1$-Hausdorff measure).

For a hypersurface $S \subseteq \mathbb{R}^k$ with its singular sets having zero hypersurface measure, we can define the \emph{directed area} function $A_S : S^{k-1} \rightarrow \mathbb{R}$ such that:
\begin{equation}\label{directedareaequation1}
A_S (v) : = \int_S |v \cdot n| \mathrm{d} Area_S
\end{equation}

Here $n$ is the normal vector that makes sense almost everywhere.

$A_S (v)$ has another formula which is also easy to understand \cite{guth2010endpoint}: Let $\pi_v : \mathbb{R}^k \rightarrow v^{\perp}$ be the orthogonal projection onto the hyperplane $v^{\perp}$. Then

\begin{equation}\label{directedareaequation2}
A_S (v) = \int_{v^{\perp}} |S \cap \pi_v^{-1} (y)| \mathrm{d} y.
\end{equation}

We now recall the cylinder estimate and a large directed area lemma which are proved in, say, \cite{guth2010endpoint}. The cylinder estimate can be easily proved using (\ref{directedareaequation2}) while the other lemma is obvious from (\ref{directedareaequation1}).

\begin{lem}[Cylinder estimate]\label{cylinderestimatelemma}
In dimension $k$, for a cylinder $T$ of radius $R$, a unit vector $v$ parallel to the axis of $T$ and an algebraic surface $S$ of degree $d$, we have
\begin{equation}
A_{S \cap T} (v) \lesssim_k R^{k-1} d.
\end{equation}
\end{lem}

\begin{lem}[A surface with large area must have large directed area for some direction]\label{largedirectedarealemma}
In dimension $k$, we let $e_1, \ldots, e_k$ be the unit coordinate vectors. Let $S$ be a hypersurface in $\mathbb{R}^k$, then
\begin{equation}
Area (S) \leq \sum_{i=1}^k A_{S} (e_i).
\end{equation}
\end{lem}

With the background knowledge above, we now proceed to the proof. First we prove a very general proposition about arbitrary hypersurfaces by ``induction on dimension''. The proposition intuitively asserts that a surface either has large area or is ``degenerate'' on some dimensions. For coordinate vectors $e_i, i \in I  \subseteq \{1, 2, \ldots, k\}$, we use the convention $< e_i : i \in I>$ to denote the coordinate subspace $\{\sum_{i \in I} t_i e_i : t_i \in \mathbb{R} \}$ generated by $\{e_i : i \in I \}$. Also, for $t > 0$ we denote $t I^n$ to be the cube of edge length $t$ with the same center as the unit cube $I^n$.

\begin{prop}\label{inductionofdimprop}
Given any $0 < t < 1$, dimension $k \geq 2$. Assume $S$ is a (closed) hypersurface (probably with boundary) contained in the $k$ dimensional unit cube $I^k$ with its singular points having zero hypersurface measure. Also assume $S \cap t I^k \neq \emptyset$. Then the implied constants can be chosen such that either $Area (S) \gtrsim_{k, t} 1$ or that there exists a coordinate subspace $W = <e_i : i \in I_W \subseteq \{1, 2, \ldots, k\}> \subseteq{R}^k$ ($I_W \neq \emptyset$), and a measurable set $B \subseteq <e_i : i \notin I_W>$ satisfying the following two conditions: (the $k-|I_W|$-dimensional measure) $|B| \gtrsim_{k, t} 1$; For an arbitrary $b \in B$, the translated hypersurface $W_b = W + b$ satisfies that there is a connected component of $W_b \cap S$ contained in the interior $I^k \setminus \partial I^k$. In the case that $W = \mathbb{R}^k$ we simply mean by the above conditions that (instead of requiring the existence of $B$) there is a connected component of $S$ contained in the interior $I^k \setminus \partial I^k$.
\end{prop}

\begin{proof}
As advertised, we use an ``induction on dimension'' argument. If the dimension $k = 2$ then $S$ is just a curve and has a point inside $t I^2$. Thus either the connected component of this point is not shorter than $\frac{1-t}{2}$ or that the component does not touch $\partial I^2$. In the latter case we take $W$ to be the whole space and the proposition holds.

Assume we already proved the proposition for $k = k_0 \geq 2$. Now assume the dimension $k  = k_0 +1$. Take any point $p \in S \cap t I^{k_0 + 1}$. If the connected component $Z_1$ of $S$ containing $p$ does not touch $\partial ((\frac{t+1}{2})I^{k_0 +1})$ then we are done.

From now on we assume that $Z_1 \cap \partial ((\frac{t+1}{2})I^{k_0 +1}) \neq \emptyset$. Then without loss of generality we may assume there is a point $p_{\frac{t+1}{2}} \in Z_1 \cap ((\frac{t+1}{2})I^{k_0 +1})$ and that the first coordinate of $p_{\frac{t+1}{2}}$ is $\frac{t+1}{2}$. Since $Z_1$ is connected, we may further assume that for any $t \leq x \leq \frac{t+1}{2}$, we can find a point $p_x \in Z_1 \cap ((\frac{t+1}{2})I^{k_0 +1})$ with the first coordinate of $p_x$ being $x$.Now apply the induction hypothesis with dimension $k_0$ with the parameter $t$ replaced by $\frac{t+1}{2}$ to the $k_0$-dimensional cube $ I^{k_0 }_x = I^{k_0 +1} \cap \{x_1 = x\} (t \leq x \leq \frac{t+1}{2})$. Note that there are only finite many cases in the conclusion of the induction hypothesis. We deduce that there are some choices of implied constants only depending on $t$ and $k_0$ such that there is a measurable subset $E \subseteq [t, \frac{t+1}{2}] \subseteq \mathbb{R}$ satisfying $|E| \gtrsim_{k_0 , t} 1$ and that for all $I^{k_0}_x, x \in E$, one same case in the conclusion of the induction hypothesis always holds.

Assume first that the $k_0-1$ dimensional area of $S \cap I^{k_0}_x$ (it is finite for all but finite $x$ and we just exclude those ``bad'' $x$ from $E$) is $\gtrsim_{k_0, t} 1$ for all $x \in E$. Then by Lemma \ref{largedirectedarealemma} we must have a $e_j (2 \leq j \leq k_0 + 1)$ and a $E' \subseteq E$, $|E'| \geq \frac{1}{k_0} |E|\gtrsim_{k_0, t} 1$ such that for each $x \in E'$, the directed ($k_0 -1$ dimensional) area $A_{S \cap I^{k_0}_x} (e_j)$ (within the hyperplane where $I^{k_0}_x$ lies) is $\gtrsim_{k_0, t} 1$. Hence for every $x \in E'$ the projection of $S\cap I^{k_0}_x$ onto the orthogonal complement of the direction $e_j$ (within $I^{k_0}_x$) has $k_0 -1$ dimensional measure $\gtrsim_{k_0, t} 1$ by (\ref{directedareaequation2}). Now that the projection of the whole $S$ onto the orthogonal complement of $e_j$ in the whole $\mathbb{R}^{k_0 + 1}$ is of course measurable and because of the above properties of $E'$, we deduce by (\ref{directedareaequation2}) and Fubini that $A_S (e_j) \gtrsim_{k_0, t} 1$ and thus $Area (S) \gtrsim_{k_0, t} 1$.

For the other cases, our arguments will be similar. First we introduce some notations. For an affine coordinate subspace $W$ (of arbitrary dimension), we say that $W$ satisfies the ``isolation condition'' (with respect to $S$ and $I^{k_0 + 1}$) if $W \cap S$ has a connected component contained in the interior of the cube slice $W \cap I^{k_0 + 1}$. By convention, we say that this condition fails if $W \cap I^{k_0 + 1} = \emptyset$. If the case in the above paragraph does not apply, then by the induction hypothesis we may assume that there is an index set $I \subseteq \{2, 3, \ldots, k_0 + 1\}$ such that for any $x \in E$ and any fixed affine subspace $W_{I, x}$ inside $I^{k_0}_x$ parallel to $<e_i : i \in I>$ (and having the same dimension as $<e_i : i \in I>$), the ($k_0 - |I|$-dimensional) measure of $Y_x = \{(\alpha_i), i \in \{2, 3, \ldots, k_0 +1\} \setminus I : \alpha_i \in \mathbb{R}, W_{I, x} + \sum_{i \in  \{2, 3, \ldots, k_0 +1\} \setminus I} \alpha_i e_i \text{ satisfies the isolation condition}\}$ is $\gtrsim_{k_0 , t} 1$. We now look at the set $Y = \{(\alpha_i), i \in \{1, 2, 3, \ldots, k_0 +1\} \setminus I : \alpha_i \in \mathbb{R}, <e_i : i \in I> + \sum_{i \in  \{1, 2, 3, \ldots, k_0 +1\} \setminus I} \alpha_i e_i \text{ satisfies the isolation condition}\}$. By quantize the isolation condition (i.e. looking at whether there is a component of $W \cap S$ that has its distance to $\partial (W \cap I^{k_0 + 1})$ not less than $\frac{1}{m}$ for $m = 2, 3,  \ldots$) we see that $Y$ (as well as the previous $Y_x$) is a countable union of closed sets and hence measurable. By Fubini, $|Y| \geq \int_{E} |Y_x| \mathrm{d} x \gtrsim_{k_0 , t} 1$. We take $W = <e_i : i \in I>$ and $B = \{\sum_{i \in  \{1, 2, 3, \ldots, k_0 +1\} \setminus I} \alpha_i e_i : (\alpha_i) \in Y\}$. The conclusion then holds.
\end{proof}

Now we are ready to prove Theorem \ref{theoremondensepolynomial}.

\begin{proof}[Proof of Theorem \ref{theoremondensepolynomial}]
Without loss of generality, we may assume $\alpha d \leq 1$ and that $\frac{1}{\alpha}$ is an integer divisible by $100$. We denote the polynomial by $Q$ and denote the hypersurface $\{Q = 0\}$ by $S$. Let $X$ be a large positive constant which will be chosen later. Our $X$ will be independent of $\alpha$ and $d$ but only depends on the dimension $n$. Let $\gamma = 10X\alpha$. All the cubes here will be closed cubes.

We divide $I^n$ into $\alpha^{-n} = (\frac{10X}{\gamma})^n$ equal cubes of edge length $\frac{\gamma}{10X}$. For each such small cube $\Delta_i$, denote $2 \Delta_i$ to be the cube of edge length $\frac{\gamma}{5X}$ with the same center as $\Delta_i$. For each $\Delta_i$, either there is a zero of $Q$ inside $2\Delta_i$ or that the $\frac{\gamma}{10X}$-neighborhood of $S$ does not intersect the interior of $\Delta_i$. Thus if we redivide $I^n$ into $(\frac{X}{\gamma})^n$ equal cubes $\Sigma_j$ of edge length $\frac{\gamma}{X}$. Then there exists a constant $s_n$ such that either that the $\frac{\gamma}{10X}$-neighborhood of $S$ intersects the interior of no more than $(\frac{10X}{\gamma})^n \gamma d$ small cubes $\Delta_i$ (and thus has a volume $\leq \gamma d \lesssim_{n} \alpha d$, since that $X$ only depends on $n$), or that $S$ intersects $\geq s_n (\frac{X}{\gamma})^n \gamma d$ cubes $\Sigma_j$. From now on, we always assume the latter case happens since in the former case we are done.

Obviously, by refining the set $\{\Sigma_j\}$ and a slight changing of $s_n$, we may assume further that $S$ intersects $\geq s_n (\frac{X}{\gamma})^n \gamma d$ cubes $\Sigma_j$ such that the cubes $2\Sigma_j$ are all contained inside $I_n$ and each pair $2\Sigma_{j_1}$ and $2\Sigma_{j_2}$ have no interior points in common. Invoke Proposition \ref{inductionofdimprop} (with $t = \frac{1}{2}$, $k = n$ and $I^n$ replaced by a rescaled and translated copy of $2 \Sigma_j$) and classify the $j$'s according to the case that applies to $2\Sigma_j$ in the conclusion of the proposition. By pigeonhole we deduce that there is an index set $\Lambda$, $|\Lambda| \gtrsim_{n} (\frac{X}{\gamma})^n \gamma d$ such that one of the following three assertions always holds for $j \in \Lambda$:

1) $Area (S \cap 2 \Sigma_j) \gtrsim_{n} (\frac{\gamma}{X})^{n-1}$ for every $j \in \Lambda$;

2) For each $2 \Sigma_j$, $j \in \Lambda$, there is a connected component of $S \cap 2\Sigma_j$ that does not touch $\partial (2 \Sigma_j)$.

3) There exists a fixed index set $\emptyset \subsetneqq I \subsetneqq \{1, 2, \ldots, n\}$ (depending on the whole $\Lambda$ but independent of $j$) such that there exists a $B_j \subseteq <e_i: i \notin I>$ for every $j \in \Lambda$ satisfying that (the $n-|I|$-dimensional measure) $|B_j| \gtrsim_{n} (\frac{\gamma}{X})^{n-|I|}$ and that for any $b \in B_j$, there is a connected component of $S \cap W_b \cap 2\Sigma_j$ that does not touch $\partial (W_b \cap 2 \Sigma_j)$, where $W_b$ is the affine subspace $<e_i: i \in I> + b$. We require $S \cap W_b \cap 2\Sigma_j \neq \emptyset$ for all $j \in \Lambda$ in this case.

If we are in case 1), then by Lemma \ref{largedirectedarealemma}, there exists an $e_i$ such that
\begin{equation}
A_{S \cap I^n} (e_i) \gtrsim_{n} Xd.
\end{equation}

But by the cylinder estimate (Lemma \ref{cylinderestimatelemma}),
\begin{equation}
A_{S \cap I^n} (e_i) \lesssim_{n} d.
\end{equation}

In fact the $\lesssim_{n}$ is a $\leq$ here but we do not need it. From here we already see that if we take $X$ large depending only on $n$ then the first case does not occur.

If we are in case 2), then since the interior of $2\Sigma_j$ do not intersect, we deduce that $S$ must have $\gtrsim_{n} (\frac{X}{\gamma})^n \gamma d$ connected components. But by the generalized Harnack curve Theorem \ref{generalizedharnack}, we know that the number of connected components is $\lesssim_n d^n$. Since $ \frac{\gamma d}{10X} = \alpha d \leq 1$, the second case does not occur if we take $X$ large depending only on $n$.

If we are in case 3), then we first note that each $B_j$ has to be contained inside the unit cube in the subspace $<e_i : i \notin I>$. This cube has $n-|I|$-dimensional measure $1$. Thus there exists a $b \in <e_i : i \notin I>$ which belongs to $\gtrsim_{n, J} (\frac{X}{\gamma})^{|I|} \gamma d$ different $B_j$ (and that $Q$ does not vanish identically on $W_b = <e_i : i \in I> + b$ since the ``bad'' $b$'s trivially have measure zero). Hence the affine subspace $W_b$ of dimension $|I|$ such that the intersection of $S$ and $W_b$ has $\gtrsim_{n, J} (\frac{X}{\gamma})^{|I|} \gamma d$ connected components. However, by the generalized Harnack again, the number of the connected components should not exceed $O_{n} (d^{|I|})$. Therefore, again with the fact that $ \frac{\gamma d}{10X} \leq 1$, we deduce that if we take $X$ large depending only on $n, J$ then the third case does not occur, either.

Hence as long as we take $X$ large depending on $n$ and $J$, all the three above cases do not happen. Thus the only case left with us is the case we excluded at the beginning that the $\alpha$-neighborhood of $S$ intersects the interior of no more than $(\frac{10X}{\gamma})^n \gamma d$ small cubes $\Delta_i$ (and thus has a volume $\leq \gamma d \lesssim_n \alpha d$). Whence the theorem holds.
\end{proof}

\section{A negative answer to Question \ref{naiveques}}

In this section, we give a negative answer to the previous Question \ref{naiveques} for all sufficiently large dimensions $n$. We will take integer points and lines lying on a quadratic hypersurface and show that they have higher incidental pattern than what would be predicted by a positive answer to Question \ref{naiveques}. In another recent work \cite{solomon2014highly}, the situation is somewhat similar (though the techniques are different).

More specifically, we will prove the following theorem.

\begin{thm}\label{counterexamplethm}
For all sufficiently large $n$ and all $N$, there is a $\delta (n) > 0$, such that there exists a point set $P \subseteq \mathbb{R}^n$ and a line set $L$, satisfying: there are $\lesssim_n N^{r-1}$ lines lying in a common $r$-dimensional linear subspace ($1 < r < n$); each line in $L$ passes through $\sim_n N$ points, $|L| \sim_n N^{n-1}$ and that $|P| \lesssim_n N^{n - \delta (n)}$.
\end{thm}

\begin{proof}
In the proof we always assume that the dimension $n$ is sufficiently large. Also we can assume that for a fixed $n$,  the parameter $N$ is large enough. We choose a nondegenerate quadratic form $Q$ over $\mathbb{R}^n$ with positive and negative inertia indices both close to $\frac{n}{2}$ (could be $\frac{n-1}{2}, \frac{n}{2}$ or $\frac{n-1}{2}$ according to the parity of $n$). For simplicity we also use $Q$ to denote the bilinear form associated with the quadratic form $Q$. The definition of $Q$ will be precise in the rest of the proof according to the number of variables it has.

Now take the set $P$ of points to be $\{x \in \mathbb{Z}^n : |x| \leq N^{1+\alpha}, Q(x) = 1\}$. Here $\alpha > 0$ is to be determined. We can use the circle method \cite{birch1962forms} to count $|P|$ when $n$ is large. Note that the ``singular locus'' defined in \cite{birch1962forms} is empty, we deduce that $|P| \sim_n N^{n + n\alpha - 2 - 2\alpha}$.

For the lines, we take $L$ to be $\{x + tv: x \in \mathbb{Z}^n, v \in \mathbb{Z}^n, |x| \leq N^{1+\alpha}, \frac{N^{\alpha}}{2} \leq |v| \leq N^{\alpha}, Q(x + t v) \equiv 1\}$. This construction is analogous to a construction of Tao \cite{tao2005new}. The last condition is equivalent to $Q(x, x) = 1$, $Q(x, v) = Q(v, v) = 0$. This is a quadratic system and we again apply the circle method. It is easy to see that any point $(x_0, v_0)$ in the singular locus satisfies $v_0 = 0$. This has dimension $n - 1$ in $\mathbb{R}^{2n}$. We claim that a circle method similar to that in \cite{birch1962forms} applies when $n$ is large, which will be explained in the next paragraph. We have immediately the fact that each line passes through $\sim_n N$ points. If the circle method works, we would deduce that $|L| \sim_n N^{n + 2n \alpha - 4 - 6\alpha}$.

We stop for a while and justify the circle method used in the last paragraph. We have to count $|\{(x, v): x, v \in \mathbb{Z}^{n}, |x| \leq N^{1+\alpha}, \frac{N^{\alpha}}{2} \leq |v| \leq N^{\alpha}, Q(x, x) = 1, Q(x, v) = Q(v, v) = 0\}|$. Cut this set into small sets $\{(x, v): x, v \in \mathbb{Z}^{n}, |x - x_j| \leq N^{\alpha}, \frac{N^{\alpha}}{2} \leq |v| \leq N^{\alpha}, Q(x, x) = 1, Q(x, v) = Q(v, v) = 0\}$ where the $x_j$'s form a grid of unit side length $N^{\alpha}$. For each small set we can count the number of points using Birch's method and results. Note that the quadratic coefficients are unchanged under translation, so the following two terms both can be uniformly controlled: (1) the upper bound of things on the minor arcs and (2) the difference between things on the major arcs and the main terms. Now look at all the main terms and we found that the difference of them lie only in the singular integral. We sum all the singular integral up, and in order that this converges to one single main term, we have to show that the following limit exists:
\begin{equation}\label{countingmaintermlimit}
\lim_{N \rightarrow \infty} \lim_{\beta \rightarrow \infty} N^{-n-n\alpha} \int_{|\gamma| \leq \beta} \int_{|x| \leq N, |v| \sim N^{\alpha}} e (\gamma_1 (Q(x, x)-1) + \gamma_2 Q(x, v) + \gamma_3 Q(v, v)) \mathrm{d} x \mathrm{d} v \mathrm{d} \gamma.
\end{equation}
Here $e(\cdot ) = e^{2 \pi \mathrm{i} \cdot}$. After scaling on the variable $v$, this integral is essentially transformed to the singular integral considered in \cite{birch1962forms}. So we know that it indeed has a limit.

We choose $\alpha$ such that $|L| \sim |N|^{n-1}$. So we have to take $\alpha = \frac{3}{2n-6}$. It suffices to prove that for any $r$-dimensional affine linear subspace ($2 \leq r \leq n-1$), there are $\lesssim N^{r-1}$ lines of $L$ lying in it. Take an arbitrary $r$-dimensional affine subspace $S_r$. Assume without loss of generality that the rank of the sublattice $\mathscr{L}_r = \mathbb{Z}^n \bigcap S_r$ is $r$. Denote the translation of $S_r$ that passes through the origin by $V_r$, and denote $\mathscr{K}_r = \mathbb{Z}^n \bigcap V_r$. Taking into account of the fact that each line in $L$ passes through $\sim N$ points in $P$, we have to prove that $|\{(x, v): x \in \mathscr{L}_r, v \in \mathscr{K}_r, |x| \lesssim N^{1 + \alpha}, |v| \sim N^{\alpha}, Q(x, x) = 1, Q(x, v) = Q(v, v) = 0\}| \lesssim N^r$. We need to be careful that this bound has to be independent of $S_r$. All the estimates below will be of this kind.

Because that the shortest distance between two points in $\mathscr{K}_r$ (or $\mathscr{L}_r$) is $\gtrsim 1$, we deduce that there is an affine linear transform that transforms $\mathscr{K}_r$ bijectively onto $\mathbb{Z}^r$, and that transforms the ball $\{v \in S_r: |v| \lesssim N^{\alpha} \}$ into the ball $\{|y| \lesssim N^{\alpha}\}$ with some carefully chosen implied constant. Theorem 2 in \cite{heath2002density} shows that $|\{v \in \mathscr{L}_r: |x| \lesssim N^{\alpha}, Q(x, x) = 1\}| \lesssim N^{(r-2) \alpha + \varepsilon}$, as long as $r \geq 3$.

For any such fixed $v$, the vectors $x \in \mathscr{L}_r$ s.t. $Q(x, v) = 0$ form a sublattice of rank $r-1$ in $\mathscr{L}_r$ (notice that $\mathscr{L}_r$ is a translated copy of $\mathscr{K}_r$). Using the same method as in the last paragraph again, we obtain that the number of vectors $x$ in this sublattice such that $Q(x, x) = 1$ is $\lesssim N^{r-3+(r-3) \alpha + \varepsilon}$, as long as $r \geq 4$. Notice that here we need a generalization of Theorem 2 in \cite{heath2002density} from quadratic forms to quadratic polynomials. It can be obtained immediately from the original proof.

In total we obtain that $|\{(x, v): x \in \mathscr{L}_r, v \in \mathscr{K}_r, |x| \lesssim N^{1 + \alpha}, |v| \sim N^{\alpha}, Q(x, x) = 1, Q(x, v) = Q(v, v) = 0\}| \lesssim N^{r-3 + (2r -5) \alpha + \varepsilon}$. This is $ < N^{r}$ since $2r-5<2n-6$ when $r \leq n-1$.

We are left with the case when $r = 2$ and $r = 3$. When $r = 3$, we can reverse the above procedure and first control the number of $x$ with $Q(x, x) = 1$. Then we simply use the trivial estimate $|\{v \in \mathscr{K}_3 : |v| \sim N^{\alpha}\}| \lesssim N^{3 \alpha}$. We deduce that $|\{(x, v): x \in \mathscr{L}_3, v \in \mathscr{K}_3, |x| \lesssim N^{1 + \alpha}, |v| \sim N^{\alpha}, Q(x, x) = 1, Q(x, v) = Q(v, v) = 0\}| \lesssim N^{1 + 4\alpha + \varepsilon}$. When $n$ is large this is definitely $\lesssim N^2 \lesssim N^3$. And we can use this bound to deal with the case $r=2$, too.
\end{proof}
\bibliographystyle{amsalpha}
\bibliography{stageref}

Department of Mathematics, Princeton University, Princeton, NJ 08540

ruixiang@math.princeton.edu

\end{document}